\documentclass[10pt, intlimits, reqno]{amsart}
\usepackage{amsmath, amssymb, enumerate, mathrsfs}

\DeclareMathOperator{\im}{Im}
\DeclareMathOperator{\sh}{sh}
\newcommand{\rme}{\mathrm{e}}

\begin{document}
\title[Trace formula for the Sturm--Liouville equation]{The regularized trace formula for the Sturm--Liouville equation with spectral parameter in the boundary conditions}
\author{Namig~J.~Guliyev}
\begin{abstract}
  The regularized trace formula of first order for the Sturm--Liouville equation with spectral parameter in the boundary conditions is obtained.
\end{abstract}
\maketitle

Consider the boundary value problem
\begin{align}
  -y''(x) + q(x)y(x) & = \lambda y(x), \qquad x \in [0,\pi], \label{eq:main} \\
  \lambda (y'(0)   - h y(0))   & = h_1 y'(0)   - h_2 y(0), \label{eq:boundary1} \\
  \lambda (y'(\pi) + H y(\pi)) & = H_1 y'(\pi) + H_2 y(\pi), \label{eq:boundary2}
\end{align}
where $q(x) \in \mathscr{W}_2^1[0,\pi]$ is a real-valued function, $h, h_1, h_2, H, H_1, H_2 \in \mathbb{R}$ and
\begin{equation*}
  \delta := h h_1 - h_2 > 0, \quad \rho := H H_1 - H_2 > 0.
\end{equation*}
Let $\varphi(x, \lambda)$ be the solution of~\eqref{eq:main} satisfying the initial conditions
\begin{equation*}
  \varphi(0, \lambda) = h_1 - \lambda, \quad \varphi'(0, \lambda) = h_2 - \lambda h.
\end{equation*}
We denote $s := \sqrt{\lambda}$. Since $\varphi(x, \lambda)$ satisfies the integral equation
\begin{equation*}
  \varphi(x, \lambda) = (h_1 - s^2) \cos sx + (h_2 - hs^2) \frac{\sin sx}{s} + \int_0^x \frac{\sin [s(x-t)]}{s} q(t) \varphi(t, \lambda) dt
\end{equation*}
and the estimate
\begin{equation*}
  \varphi(x, \lambda) = O \left( |s|^2 e^{|\im s|x} \right),
\end{equation*}
one can obtain the following representations:
\begin{equation} \label{eq:phi}
  \varphi(x, \lambda) = -s^2 \cos sx - \left( h + \frac{1}{2} \int_0^x q(t)dt \right) s \sin sx + O \left( e^{|\im s|x} \right),
\end{equation}
\begin{multline} \label{eq:phi'}
  \varphi'(x, \lambda) = s^3 \sin sx - \left( h + \frac{1}{2} \int_0^x q(t)dt \right) s^2 \cos sx - \\
  - \left( h_1 + \frac{q(x) + q(0)}{4} + \frac{h}{2} \int_0^x q(t)dt + \frac{1}{8} \left( \int_0^x q(t)dt \right)^2 \right) s \sin sx - \\
  - \frac{s}{4} \int_0^x q'(t) \sin [s(x-2t)] dt + O \left( e^{|\im s|x} \right).
\end{multline}
Denote by $\chi(\lambda)$ the characteristic function of~\eqref{eq:main}--\eqref{eq:boundary2}:
\begin{equation*}
  \chi(\lambda) := (\lambda - H_1) \varphi'(\pi, \lambda) + (\lambda H - H_2) \varphi(\pi, \lambda).
\end{equation*}
From~\eqref{eq:phi} and~\eqref{eq:phi'} we have
\begin{multline} \label{eq:chi1}
  \chi(\lambda) := s^5 \sin s\pi - \omega s^4 \cos s\pi - \omega_1 s^3 \sin s\pi - \frac{s^3}{4} \int_0^{\pi} q'(t) \sin [s(\pi-2t)] dt + \\
  + O \left( |s|^2 \rme^{|\im s|\pi} \right),
\end{multline}
where
\begin{equation*}
  \omega := h + H + \frac{1}{2} \int_0^{\pi} q(t)dt,
\end{equation*}
\begin{equation*}
  \omega_1 := h_1 + H_1 + hH + \frac{q(\pi) + q(0)}{4} + \frac{h+H}{2} \int_0^{\pi} q(t)dt + \frac{1}{8} \left( \int_0^{\pi} q(t)dt \right)^2.
\end{equation*}
The eigenvalues of the problem~\eqref{eq:main}--\eqref{eq:boundary2} are real and coincide with the zeros of the function $\chi(\lambda)$. From~\eqref{eq:chi1} we obtain
\begin{equation*}
  s_n = n - 2 + \frac{\omega}{\pi(n-2)} + \frac{\zeta_n}{n^2}, \quad \{\zeta_n\} \in l_2.
\end{equation*}
Then
\begin{equation} \label{eq:lambda_n}
  \lambda_n = s_n^2 = (n - 2)^2 + \frac{2\omega}{\pi} + \frac{\xi_n}{n}, \quad \{\xi_n\} \in l_2.
\end{equation}
Therefore
\begin{equation} \label{eq:s_lambda}
  s_{\lambda} := \lambda_0 + \lambda_1 + \sum_{n=2}^{\infty} \left( \lambda_n - (n-2)^2 - \frac{2\omega}{\pi} \right) < \infty.
\end{equation}
The series~\eqref{eq:s_lambda} is called {\it the regularized trace of first order} (or {\it the first regularized trace}) for the problem~\eqref{eq:main}--\eqref{eq:boundary2}. Analogously, one can define {\it the second}, {\it the third} etc. regularized traces. Note that the second regularized trace formula for the Sturm--Liouville problem with spectral parameter in a boundary condition has been found in \cite{Albayrak2000}.

Since $\chi(\lambda)$ is an entire function of order $1/2$, from Hadamard's theorem (see \cite[Section~4.2]{Levin1996}), using~\eqref{eq:chi1} we obtain
\begin{equation} \label{eq:chi2}
  \chi(\lambda) = \pi(\lambda-\lambda_0)(\lambda-\lambda_1)(\lambda-\lambda_2) \prod_{n=3}^{\infty} \frac{\lambda_n - \lambda}{(n-2)^2}
\end{equation}

Let $\lambda= -\mu^2$. We calculate the sum $s_{\lambda}$ of the series~\eqref{eq:s_lambda} by comparing the asymptotic expressions obtained from formulas~\eqref{eq:chi1} and~\eqref{eq:chi2} as $\mu \to \infty$.
From the formula~\eqref{eq:chi2} we have
\begin{multline} \label{eq:chi3}
  \chi(-\mu^2) = -\pi(\mu^2+\lambda_0)(\mu^2+\lambda_1)(\mu^2+\lambda_2) \prod_{n=3}^{\infty} \frac{\lambda_n + \mu^2}{(n-2)^2} = \\
  = - \frac{(\mu^2+\lambda_0)(\mu^2+\lambda_1)(\mu^2+\lambda_2) \sh \pi\mu}{\mu}\Phi(\mu),
\end{multline}
where
\begin{equation*}
  \Phi(\mu) = \prod_{n=1}^{\infty} \left( 1 - \frac{n^2 - \lambda_{n+2}}{\mu^2 + n^2} \right).
\end{equation*}
We need the formulas (see~\cite{Levitan1964})
\begin{equation} \label{eq:Levitan1}
  \sum_{k=2}^{\infty} \frac{1}{k} \sum_{n=1}^{\infty} \frac{|n^2 - \lambda_{n+2}|^k}{(\mu^2 + n^2)^k} = o\left( \frac{1}{\mu^2} \right)
\end{equation}
and
\begin{equation} \label{eq:Levitan2}
  \sum_{n=1}^{\infty} \frac{1}{\mu^2 + n^2} = \frac{\pi}{2\mu} - \frac{1}{2\mu^2} + o\left( \frac{1}{\mu^2} \right).
\end{equation}
Using~\eqref{eq:lambda_n} and~\eqref{eq:Levitan2} we have
\begin{multline} \label{eq:temp}
  \frac{1}{\mu^2} \sum_{n=1}^{\infty} \left( \lambda_{n+2} - n^2 - \frac{2\omega}{\pi} \right) \frac{n^2}{\mu^2 + n^2} = \frac{1}{\mu^2} \sum_{n=1}^{\infty} \frac{n\xi_{n+2}}{\mu^2 + n^2} \le \\
  \le \frac{1}{\mu^2} \sqrt{\sum_{n=3}^{\infty} \xi_n^2} \sqrt{\sum_{n=1}^{\infty} \frac{n^2}{(\mu^2 + n^2)^2}} \le \frac{\Omega}{\mu^2} \sqrt{\sum_{n=1}^{\infty} \frac{\mu^2 + n^2}{(\mu^2 + n^2)^2}} = \\
  = \frac{\Omega}{\mu^2} \sqrt{\sum_{n=1}^{\infty} \frac{1}{\mu^2 + n^2}} = o\left( \frac{1}{\mu^2} \right),
\end{multline}
where
\begin{equation*}
  \Omega := \sqrt{\sum_{n=3}^{\infty} \xi_n^2} < \infty.
\end{equation*}
From~\eqref{eq:Levitan1}, \eqref{eq:Levitan2} and~\eqref{eq:temp} we calculate:
\begin{multline*}
  \ln \Phi(\mu) = \sum_{n=1}^{\infty} \left( 1 - \frac{n^2 - \lambda_{n+2}}{\mu^2 + n^2} \right) = - \sum_{n=1}^{\infty} \sum_{k=1}^{\infty} \frac{1}{k} \left( \frac{n^2 - \lambda_{n+2}}{\mu^2 + n^2} \right)^k = \\
  = - \sum_{k=1}^{\infty} \frac{1}{k} \sum_{n=1}^{\infty} \left( \frac{n^2 - \lambda_{n+2}}{\mu^2 + n^2} \right)^k = - \sum_{n=1}^{\infty} \frac{n^2 - \lambda_{n+2}}{\mu^2 + n^2} + o\left( \frac{1}{\mu^2} \right) = \\
  = \frac{2\omega}{\pi} \sum_{n=1}^{\infty} \frac{1}{\mu^2 + n^2} + \sum_{n=1}^{\infty} \left( \lambda_{n+2} - n^2 - \frac{2\omega}{\pi} \right) \frac{1}{\mu^2 + n^2} + o\left( \frac{1}{\mu^2} \right) = \frac{\omega}{\mu} - \frac{\omega}{\pi\mu^2} + \\
  + \frac{1}{\mu^2} \sum_{n=1}^{\infty} \left( \lambda_{n+2} - n^2 - \frac{2\omega}{\pi} \right) - \frac{1}{\mu^2} \sum_{n=1}^{\infty} \left( \lambda_{n+2} - n^2 - \frac{2\omega}{\pi} \right) \frac{n^2}{\mu^2 + n^2} + \\
  + o\left( \frac{1}{\mu^2} \right) = \frac{\omega}{\mu} - \frac{\omega}{\pi\mu^2} + \frac{1}{\mu^2} \left( s_{\lambda} - \lambda_0 - \lambda_1 - \lambda_2 + \frac{2\omega}{\pi} \right) + o\left( \frac{1}{\mu^2} \right) = \\
  = \frac{\omega}{\mu} + \frac{1}{\mu^2} \left( s_{\lambda} - \lambda_0 - \lambda_1 - \lambda_2 + \frac{\omega}{\pi} \right) + o\left( \frac{1}{\mu^2} \right).
\end{multline*}
Therefore
\begin{equation*}
  \Phi(\mu) = 1 + \frac{\omega}{\mu} + \frac{1}{\mu^2} \left( s_{\lambda} - \lambda_0 - \lambda_1 - \lambda_2 + \frac{\omega}{\pi} + \frac{\omega^2}{2} \right) + o\left( \frac{1}{\mu^2} \right)
\end{equation*}
and from~\eqref{eq:chi3}
\begin{equation} \label{eq:chi4}
  \chi(-\mu^2) = \frac{\rme^{\pi\mu}}{2} \left( -\mu^5 - \omega\mu^4 + \mu^3 \left( s_{\lambda} + \frac{\omega}{\pi} + \frac{\omega^2}{2} \right) + o \left( \mu^3 \right) \right).
\end{equation}
The estimate~\eqref{eq:chi1} implies
\begin{equation} \label{eq:chi5}
  \chi(-\mu^2) = \frac{\rme^{\pi\mu}}{2} \left( -\mu^5 - \omega \mu^4 + \omega_1 \mu^3 + o \left( \mu^3 \right) \right).
\end{equation}
Comparing~\eqref{eq:chi4} and~\eqref{eq:chi5} we obtain the following trace formula:
\begin{multline*}
  \lambda_0 + \lambda_1 + \sum_{n=2}^{\infty} \left( \lambda_n - (n-2)^2 - \frac{2\omega}{\pi} \right) = \\
  = h_1 + H_1 +  \frac{q(\pi) + q(0)}{4} - \frac{h+H}{\pi} - \frac{h^2 + H^2}{2} - \frac{1}{2\pi} \int_0^{\pi} q(t)dt.
\end{multline*}

{\bf Namig J. Guliyev}

Institute of Mathematics and Mechanics of NAS of Azerbaijan

9 F.Agayev str., AZ1141, Baku, Azerbaijan.

E-mail: njguliyev@yahoo.com

\end{document}